\documentclass[11pt]{article}

\usepackage{fullpage}
\usepackage{authblk}
\usepackage{natbib}
\usepackage{algorithm}
\usepackage{algpseudocode}
\usepackage{amsmath,amsthm,amsfonts}
\usepackage{amsmath}
\usepackage[subnum]{cases}
\usepackage{tcolorbox}
\usepackage{mathrsfs}
\usepackage{amssymb}
\usepackage{color}
\usepackage{mathrsfs}
\usepackage{empheq}
\usepackage{enumitem}
\usepackage{bm}
\usepackage{multirow}
\usepackage{booktabs}
\usepackage{makecell}
\usepackage{graphicx}
\usepackage{subcaption}
\usepackage{comment}
\usepackage{cases}
\usepackage{appendix}
\usepackage{tikz}
\usetikzlibrary{arrows,shapes}
\usepackage[colorlinks= true, linkcolor=red, citecolor=blue, urlcolor=black]{hyperref}
\usepackage{cleveref}
\usepackage{marginnote}
\usepackage{diagbox} 

\numberwithin{equation}{section}

\theoremstyle{definition}

\newtheorem{theorem}{Theorem}[section]
\newtheorem{lemma}[theorem]{Lemma}

\newtheorem{defn}[theorem]{Definition}

\newcommand{\mr}{\mathbb{R}}

\usepackage{booktabs}
\usepackage{tikz}
\usetikzlibrary{arrows.meta, positioning}
\usepackage{pgfplots}
\usetikzlibrary{calc,intersections}

\usetikzlibrary{shapes.geometric, arrows, positioning}

\tikzset{
	mybox/.style  = {draw, rectangle, minimum width=4cm, minimum height=0.8cm, text centered, text width=4.4cm,   
		font=\normalsize},
	box/.style  = {draw, rectangle, minimum width=2.0cm, minimum height=0.6cm, text centered, text width=3.0cm,   
		font=\normalsize},
	myarrow/.style = {line width=0.2pt, draw=black, -triangle 60, postaction={draw, line width=0.2pt, shorten >=10pt,-}}
}

\tikzstyle{arrow} = [->, >=stealth, -triangle 60]

\allowdisplaybreaks

\makeatletter
\newcommand{\leqnomode}{\tagsleft@true}
\newcommand{\reqnomode}{\tagsleft@false}
\makeatother

\begin{document}

\title{Lyapunov Analysis For Monotonically Forward-Backward Accelerated Algorithms}
\author[1]{Mingwei Fu}
\author[2,3]{Bin Shi\thanks{Corresponding author: \url{binshi@fudan.edu.cn} } }
\affil[1]{School of Mathematical Sciences, University of Chinese Academy of Sciences, Beijing 100049, China}
\affil[2]{Center for Mathematics and Interdisciplinary Sciences, Fudan University, Shanghai 200433, China}
\affil[3]{Shanghai Institute for Mathematics and Interdisciplinary Sciences, Shanghai 200433, China}
\date\today

\maketitle

\begin{abstract}
Nesterov’s accelerated gradient method (\texttt{NAG}) achieves faster convergence than gradient descent for convex optimization but lacks monotonicity  in function values. To address this, Beck and Teboulle [2009b] proposed a monotonic variant,~\texttt{M-NAG}, and extended it to the proximal setting as~\texttt{M-FISTA} for composite problems such as~\texttt{Lasso}. However, establishing the linear convergence of~\texttt{M-NAG} and~\texttt{M-FISTA} under strong convexity remains an open problem. In this paper, we analyze~\texttt{M-NAG} via the implicit-velocity phase representation and show that an additional assumption, either the position update or the phase-coupling relation, is necessary to fully recover the~\texttt{NAG} iterates. The essence of~\texttt{M-NAG} lies in controlling an auxiliary sequence to enforce non-increase.  We further demonstrate that the~\texttt{M-NAG} update alone is sufficient to construct a Lyapunov function guaranteeing linear convergence, without relying on full \texttt{NAG} iterates. By modifying the mixed sequence to incorporate forward-indexed gradients, we develop a new Lyapunov function that removes the kinetic energy term, enabling a direct extension to \texttt{M-NAG}. The required starting index depends only on the momentum parameter and not on problem constants. Finally, leveraging newly developed proximal inequalities, we extend our results to \texttt{M-FISTA}, establishing its linear convergence and deepening the theoretical understanding of monotonic accelerated methods.

\end{abstract}

%

\section{Introduction}
\label{sec: intro}

Over the past few decades, machine learning has emerged as a major application area for optimization algorithms. A central challenge in this field is unconstrained optimization, which involves minimizing a convex objective function without any constraints. Mathematically, this problem is formulated as:
\[
\min_{x \in \mathbb{R}^d} f(x).
\]
At the heart of addressing this challenge are gradient-based methods, which have driven significant progress due to their computational efficiency and low memory requirements. These properties make them particularly well-suited for large-scale machine learning tasks. As a result, gradient-based algorithms have become indispensable tools in modern optimization, playing a central role in both theoretical developments and real-world applications.

\begin{figure}[htb!]
\centering
\includegraphics[scale=0.24]{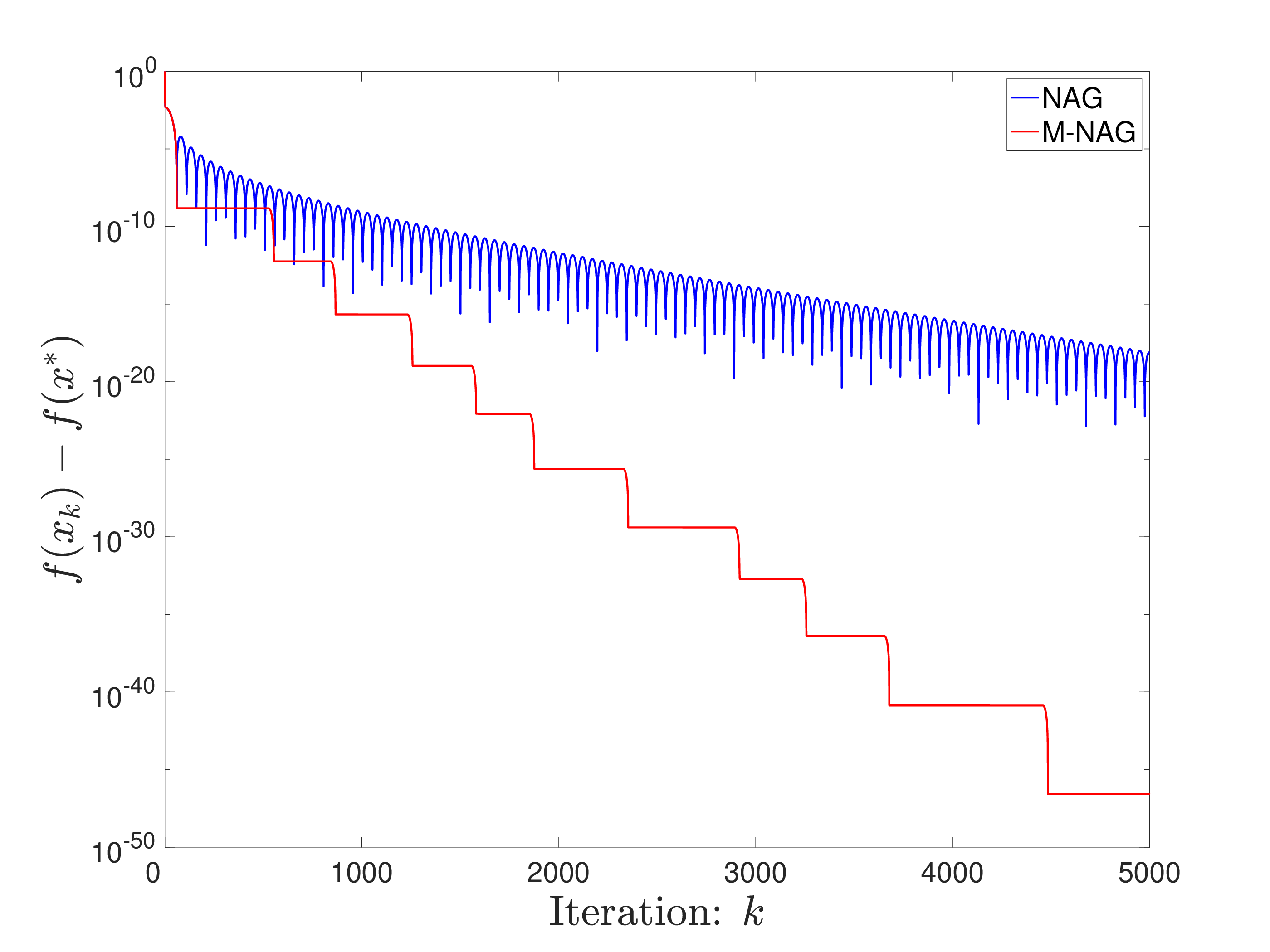}
\caption{Numerical comparison of the iterative progression of function values for~\texttt{NAG} and~\texttt{M-NAG}. The momentum parameter is set to $r=2$, and the step size is set to $s=0.4$. Both optimization algorithms are applied to the quadratic function $f(x_1, x_2) = 5 \times 10^{-3}x_1^2 + x_2^2$.}
\label{fig: quadratic}
\end{figure}

One of the most fundamental gradient-based methods is vanilla gradient descent, which dates back to~\citet{cauchy1847methode}. Given an initial point $x_0 \in \mathbb{R}^d$,  the method generates iterates according to the update rule:
\[
x_{k+1} = x_{k} - s \nabla f(x_{k}),
\]
where $s>0$ is the step size.  While this method is simple and widely used, its convergence can be slow, particularly for convex functions. To overcome this limitation,~\citet{nesterov1983method} introduced a two-step accelerated scheme, now known as Nesterov’s accelerated gradient method (\texttt{NAG}). Starting from any initial point $x_0 = y_0 \in \mathbb{R}^d$, the updates are given by:
\begin{subequations}
\label{eqn: nag}
\begin{empheq}[left=\empheqlbrace]{align}
& x_{k+1} = y_{k} - s\nabla f(y_{k}),                                                \label{eqn: nag-gradient}          \\
& y_{k+1} = x_{k+1} + \frac{k}{k+r+1} (x_{k+1} - x_{k}),                 \label{eqn: nag-momentum} 
\end{empheq}    
\end{subequations}
where $s>0$ is the step size, and $r \geq 2$ is the momentum parameter. Compared to vanilla gradient descent,~\texttt{NAG} achieves an accelerated convergence rate. However, it does not guarantee monotonicity, meaning that the function values $f(x_k)$ may not decrease at every iteration. As a result, oscillations or overshooting can occur, particularly as the iterates approach the minimizer (see~\Cref{fig: quadratic}).
This lack of stability and predictability can make it difficult to monitor optimization progress and ensure reliable convergence.

To address this issue,~\citet{beck2009fastieee} proposed a variant that ensures both acceleration and monotonicity, referred to as the Monotone Nesterov’s accelerated gradient~(\texttt{M-NAG}) method. This algorithm incorporates a comparison step that stabilizes the iterates while preserving the acceleration characteristics of~\texttt{NAG}. Starting  from any initial point $x_0 = y_0 \in \mathbb{R}^d$, the update rules of~\texttt{M-NAG} are defined as follows:  
\begin{subequations}
\label{eqn: m-nag}
\begin{empheq}[left=\empheqlbrace]{align}
& z_{k}       = y_{k} - s\nabla f(y_{k}),                                                                                                                              \label{eqn: m-nag-gradient} \\
& x_{k+1}   = \left\{ \begin{aligned} 
                               & z_{k}, && \text{if}\; f(z_{k}) \leq f(x_{k}), \\
                               & x_{k}, && \text{if}\; f(z_{k}) > f(x_{k}),
                               \end{aligned} \right.                                                                                                                          \label{eqn: m-nag-monetone} \\
& y_{k+1}  = x_{k+1} + \frac{k}{k+r+1} (x_{k+1} - x_{k}) + \frac{k+r}{k+r+1} (z_{k} - x_{k+1}),                                        \label{eqn: m-nag-momentum} 
\end{empheq}    
\end{subequations}
where $s>0$ is the step size and $r \geq 2$ is the momentum parameter. Unlike~\texttt{NAG},~\texttt{M-NAG} guarantees that the sequence of objective values $\{f(x_k)\}_{k=0}^{\infty}$ is non-increasing, leading to a more stable and predictable convergence path (see~\Cref{fig: quadratic}).  Numerical experiments by~\citet{o2015adaptive} further stimulated interest in monotonic acceleration methods, especially for strongly convex objective functions. In addition,~\citet{giselsson2014monotonicity} observed that~\texttt{M-NAG} often requires fewer iterations than~\texttt{NAG} to achieve a comparable level of accuracy.  As illustrated in~\Cref{fig: quadratic},~\texttt{M-NAG} effectively suppresses the oscillatory behavior seen in~\texttt{NAG} while preserving a faster convergence rate, in fact, it often achieves even faster convergence. These observations naturally lead to the following question:


\begin{tcolorbox}
\begin{itemize}
\item How does~\texttt{M-NAG} behave in terms of convergence when applied to strongly convex objective functions? 
\end{itemize}
\end{tcolorbox}

\subsection{Two key observations}
\label{subsec: key-observation}

In this paper, we present two key observations regarding~\texttt{NAG} and~\texttt{M-NAG}, which enable us to analyze their convergence behavior via Lyapunov techniques.
\subsubsection{Monotonicity implicitly embedded in NAG}
\label{subsubsec: monotonicity}

By substituting~\eqref{eqn: m-nag-gradient} into~\eqref{eqn: m-nag-momentum}, we find that the iterative sequences of~\texttt{M-NAG}, $\{x_k\}_{k=0}^{\infty}$ and $\{y_{k}\}_{k=0}^{\infty}$, satisfies the following relation: 
\begin{equation}
\label{eqn: essential-iteration}
y_{k+1}  = x_{k+1} +  \frac{k}{k+r+1}(x_{k+1} - x_{k}) + \frac{k+r}{k+r+1}(y_k - s \nabla f(y_k) - x_{k+1}),
\end{equation}
which shows that~\texttt{M-NAG} performs a single update  that combines information from both sequences $\{x_k\}_{k=0}^{\infty}$ and  $\{y_k\}_{k=0}^{\infty}$. More importantly, it suggests that the update of $\{y_k\}_{k=0}^{\infty}$ can be designed to guarantee that the objective values $\{f(x_k)\}_{k=0}^{\infty}$ form a non-increasing sequence. Moreover,~\eqref{eqn: essential-iteration} reveals that~\texttt{M-NAG} is essentially a linear combination of the momentum step~\eqref{eqn: nag-momentum} and the gradient step~\eqref{eqn: nag-gradient}, as depicted in~\Cref{fig: digram-nag-mnag}. In this sense,~\texttt{M-NAG} retains only part information from the full update~\texttt{NAG}~\eqref{eqn: nag}, effectively isolating and recombining its key components.  
%
%
%
%
%
%

\begin{figure}[htbp!]
\centering
\begin{tikzpicture}[
  node distance=1.2cm and 1.5cm,
  box/.style={draw, minimum width=2.4cm, minimum height=1.0cm, align=center},
  smallbox/.style={draw, minimum width=4.2cm, minimum height=1.8cm, align=center},
  every node/.style={font=\sffamily}
]

\node[smallbox, transform shape, scale=0.8] (momentum) {
  Momentum Update~\eqref{eqn: nag-momentum} \\[2pt] 
  $ y_{k+1} = x_{k+1} + \dfrac{k}{k+r+1} (x_{k+1} - x_{k}) $
};

\node[smallbox, below=3.6cm  of momentum, transform shape, scale=0.8] (gradient) {
  Gradient Update~\eqref{eqn: nag-gradient} \\[2pt] 
  $x_{k+1} = y_{k} - s \nabla f(y_k)$
};

\node[smallbox, below=0.8cm of gradient, transform shape, scale=0.8] (mnag) {\texttt{M-NAG}~\eqref{eqn: essential-iteration}\\[2pt] $y_{k+1}  = x_{k+1} +  \frac{k}{k+r+1}(x_{k+1} - x_{k}) + \frac{k+r}{k+r+1}(y_k - s \nabla f(y_k) - x_{k+1})$};
\node[box, right =2.8cm of momentum, yshift=-2.6cm, transform shape, scale=0.8] (nag) {\texttt{NAG}~\eqref{eqn: nag}};

\node[draw, circle, minimum size=1.2cm, transform shape, scale=0.8] (fraction) at ($(momentum)!0.55!(gradient) $) {\( \pmb{\dfrac{k + r}{k + r + 1}}\)};
\node at ($(fraction.north)!0.4!(momentum.south)$) {\Large$+$};  
\node at ($(fraction.south)!0.5!(gradient.north)$) {\Large$\times$};  
\node at ($(fraction.north)!1.08!(gradient.south)$) {\rotatebox{90}{\Large$=$}};

\coordinate (momentum-east) at ([xshift=3.75cm]momentum.east);
\draw[thick, -{Stealth}] (nag)-- (momentum-east) -- (momentum) ;

\coordinate (gradient-east) at ([xshift=4.65cm]gradient.east);
\draw[thick, -{Stealth}] (nag) -- (gradient-east) -- (gradient);

\end{tikzpicture}
\caption{Illustration of how~\texttt{M-NAG} consolidates the momentum and gradient steps of classical~\texttt{NAG} into a unified update rule. Each solid arrow $A \rightarrow B$ indicates that $B$ is a component or derived step of $A$.}
\label{fig: digram-nag-mnag}
\end{figure}

The acceleration mechanism of~\texttt{NAG} was rigorously uncovered by the high-resolution differential equation framework proposed by~\citet{shi2022understanding}. In particular, the implicit-velocity scheme has been shown to outperform the gradient-correction scheme in capturing sharp  convergence rates~\citep{chen2022gradient, chen2022revisiting}.  Accordingly, we adopt the phase-space representation induced by the implicit-velocity formulation. To this end, we introduce a velocity sequence $\{v_{k}\}_{k=0}^{\infty}$ defined as: 
\[
v_k = \frac{x_{k} - x_{k-1}}{\sqrt{s}}.
\] 
Under this formulation,the classical~\texttt{NAG} iteration can be equivalently rewritten as:
\begin{subequations}
\label{eqn: nag-phase}
\begin{empheq}[left=\empheqlbrace]{align}
         & x_{k+1} - x_{k} = \sqrt{s} v_{k+1},                                                                               \label{eqn: nag-phase-x}          \\
         & v_{k+1} - v_{k} = -\frac{r+1}{k+r}  v_{k} - \sqrt{s} \nabla f\left( y_{k} \right),                 \label{eqn: nag-phase-v}
\end{empheq}    
\end{subequations}
where the auxiliary sequence $\{y_{k}\}_{k=0}^{\infty}$ satisfies the phase-coupling relation:
\begin{equation}
\label{eqn: nag-phase-x-y-v}
y_k = x_{k} +  \frac{k-1}{k+r} \cdot \sqrt{s}v_k.
\end{equation}
As a key construction in their analysis,~\citet{chen2022gradient} introduced the mixed sequence
\begin{equation}
\label{eqn: mix-previous}
\mathbf{R}_k = (k-1)\sqrt{s}v_{k} + rx_k,
\end{equation}
and showed that combining~\eqref{eqn: nag-phase-x} with~\eqref{eqn: nag-phase-v} yields the update rule:
\begin{equation}
\mathbf{R}_{k+1} - \mathbf{R}_{k} = - \left( k+r\right)s \nabla f(y_{k})   \label{eqn: nag-iv-iter}. 
\end{equation}
Based on the implicit-velocity phase-space formulation given by~\eqref{eqn: nag-phase}--\eqref{eqn: nag-phase-x-y-v}, we can now make precise the connection between~\texttt{M-NAG} and~\texttt{NAG}.  
\begin{figure}[htbp!]
\centering
\begin{tikzpicture}[
  node distance=1.2cm and 1.5cm,
  box/.style={draw, minimum width=2.8cm, minimum height=1.0cm, align=center},
  smallbox/.style={draw, minimum width=4.6cm, minimum height=1.8cm, align=center},
  every node/.style={font=\sffamily},
  scale=0.68,
  transform shape
]

\node[box] (nag) {
  \texttt{NAG}~\eqref{eqn: nag}
};

\node[smallbox, below=3.6cm of nag] (nag-phase) {
 Implicit-Velocity  \\[2pt] Phase-Space Representation \\[2pt] \eqref{eqn: nag-phase} and~\eqref{eqn: nag-phase-x-y-v} 
};

\node[smallbox, right=1.6cm of nag-phase, yshift=-3.8cm] (position) {
 Position Update~\eqref{eqn: nag-phase-x} \\[2pt]
  $ x_{k+1} - x_k = \sqrt{s} v_{k+1}  $
};

\node[smallbox, below=3.8cm of position] (coupling) {
  Phase-Coupling Relation~\eqref{eqn: nag-phase-x-y-v} \\[2pt]
  $ y_{k} = x_{k} + \frac{k-1}{k+r} \cdot \sqrt{s} v_{k} $
};

\node[smallbox, right=1.6cm of nag-phase, xshift=8cm] (velocity) {
   Velocity Update~\eqref{eqn: nag-phase-v} \\[2pt]
  $ v_{k+1} - v_{k} = - \dfrac{r+1}{k+r} v_{k} - \sqrt{s} \nabla f(y_k)$
};

\node[smallbox, below=4.8cm of velocity] (mix) {
Iteration of Mixed-Sequence $\mathbf{R}_k$~\eqref{eqn: nag-iv-iter}\\[2pt] 
 $\mathbf{R}_{k+1} - \mathbf{R}_{k} = - \left( k+r\right)s \nabla f(y_{k})  $
};

\node[smallbox, below=3.8cm of mix] (mnag) {
 \texttt{M-NAG}~\eqref{eqn: essential-iteration} \\[2pt] 
 $y_{k+1}  = x_{k+1} +  \frac{k}{k+r+1}(x_{k+1} - x_{k}) + \frac{k+r}{k+r+1}(y_k - s \nabla f(y_k) - x_{k+1})$
};

Arrows
\draw[thick, - {Stealth}] (nag) --  (nag-phase);
\draw[thick, -{Stealth}] (nag-phase) -- (nag);
\draw[thick, - {Stealth}] (velocity) --  (mix);
\draw[thick, -{Stealth}] (mix) -- (velocity);
\draw[thick, - {Stealth}] (mnag) --  (mix);
\draw[thick, -{Stealth}] (mix) -- (mnag);

\draw[thick, -{Stealth}] (nag-phase) -- (velocity);

\coordinate (nag-to-pos-x) at ([yshift=-2.9cm]nag-phase.south);
\draw[thick, -{Stealth}] (nag-phase.south) -- (nag-to-pos-x) -- (position.west);

\coordinate (nag-to-coupling-x) at ([yshift=-8.5cm]nag-phase.south);
\draw[thick, -{Stealth}] (nag-phase.south) -- (nag-to-coupling-x) -- (coupling.west);

\coordinate (position-to-vm-x) at ([xshift=6.5cm]position.east);
\draw[thick, -, dashed] (position.east) -- node[midway, above]{Under Condition}(position-to-vm-x) ;

\coordinate (coupling-to-mm-x) at ([xshift=6.5cm]coupling.east);
\draw[thick, -, dashed] (coupling.east) -- node[midway, above]{Under Condition}(coupling-to-mm-x);

\end{tikzpicture}
\caption{Illustration of how~\texttt{M-NAG} recovers the~\texttt{NAG} update though the implicit-velocity phase-representation. Each solid arrow $A \rightarrow B$ indicates that $B$ is a component or derived step of $A$, while each dashed arrow represents a conditional relation required to establish the connection. }
\label{fig: digram-nag-phase-mnag}
\end{figure}
As illustrated in~\Cref{fig: digram-nag-phase-mnag}, this relationship becomes particularly transparent when examining the logical equivalence among the velocity update~\eqref{eqn: nag-phase-v}, the mixed-sequence iteration~\eqref{eqn: nag-iv-iter}, and the~\texttt{M-NAG} scheme~\eqref{eqn: essential-iteration}. Specifically:
\begin{itemize}
\item Given the position update~\eqref{eqn: nag-phase-x}, the velocity update~\eqref{eqn: nag-phase-v} is equivalent to the mixed-sequence iteration~\eqref{eqn: nag-iv-iter}. 
\item Given the phase-coupling relation~\eqref{eqn: nag-phase-x-y-v}, the mixed-sequence iteration~\eqref{eqn: nag-iv-iter} is equivalent to the~\texttt{M-NAG} update rule~\eqref{eqn: essential-iteration}. 
\end{itemize}
Since the velocity sequence $\{v_k\}_{k=0}^{\infty}$ is explicitly constructed from the position sequences $\{x_k\}_{k=0}^{\infty}$ and $\{y_k\}_{k=0}^{\infty}$, each of the two conditions, ~\eqref{eqn: nag-phase-x} and~\eqref{eqn: nag-phase-x-y-v}, can be imposed independently, without any external assumption. In this sense, both equivalences above arise naturally from the mixed-sequence iteration~\eqref{eqn: nag-iv-iter}.  However, requiring both conditions simultaneously does not follow intrinsically from \texttt{M-NAG} and must instead be imposed externally. In other words, although the \texttt{M-NAG} update~\eqref{eqn: essential-iteration} can recover the velocity update~\eqref{eqn: nag-phase-v}, and thus reconstruct the full \texttt{NAG} iteration, it requires an additional assumption: either the position update~\eqref{eqn: nag-phase-x} or the phase-coupling relation~\eqref{eqn: nag-phase-x-y-v}. This highlights that \texttt{M-NAG} only partially encodes the full structure of \texttt{NAG} derived from the implicit-velocity formulation.


\subsubsection{Lyapunov analysis and gradient iteration}
\label{subsubsec: gir}

As shown by~\citet{chen2022revisiting}, a principled approach to constructing a Lyapunov function for deriving convergence rates hinges on how the mixed sequence is defined. Since the mixed sequence is logically equivalent to the \texttt{M-NAG} update rule (as illustrated in~\Cref{subsubsec: monotonicity}), we make the important observation that the \texttt{M-NAG} scheme~\eqref{eqn: essential-iteration} forms a fundamental basis for Lyapunov analysis. Notably, this construction requires only partial information from the full \texttt{NAG} method, demonstrating that \texttt{M-NAG} alone suffices to capture the key structure needed to establish convergence.

In previous works, the mixed sequence $\{\mathbf{R}_k\}_{k=0}^{\infty}$ defined in~\eqref{eqn: mix-previous} was commonly used. For instance, in~\citet{su2016differential},~\citet{shi2022understanding}, and~\citet{chen2022gradient}, the Lyapunov function for convex objectives is typically constructed as
\begin{equation}
\label{eqn: lyapunov-cvx}
\mathcal{E}(k) = sk(k+r)\left( f(x_k) - f(x^{\star}) \right) + \frac{1}{2} \left\| (k-1)\sqrt{s}v_{k} + r(x_k - x^{\star}) \right\|^2.  \footnote{Throughout the paper, $x^{\star}$ denotes the unique global minimizer of the objective function.}
\end{equation}
If we use this Lyapunov function to derive convergence rates for strongly convex functions, one can apply the fundamental inequality of strong convexity, leading to the bound:
\begin{equation}
\label{eqn: lyapunov-cvx-iter}
\mathcal{E}(k+1) - \mathcal{E}(k) \leq - \mu s \left[ C_1(k+1)\left( f(x_{k+1}) - f(x^{\star}) \right) + C_2(k) s\left\| v_{k}  \right\|^2 + C_{3}(k)\left\| x_k - x^{\star} \right\|^2 \right],
\end{equation}
where $C_1(k) = \Theta(k^2)$, $C_2(k) = \Theta(k^2)$, and $C_3(k) = \Theta(k)$.\footnote{The fundamental inequality for strongly convex functions is rigorously stated in~\eqref{eqn: fund-inq-smooth}.}  However, comparing~\eqref{eqn: lyapunov-cvx} and~\eqref{eqn: lyapunov-cvx-iter}, we find that the right-hand side of the inequality~\eqref{eqn: lyapunov-cvx-iter} cannot be aligned proportionally with the original Lyapunov function~\eqref{eqn: lyapunov-cvx}, due to the mismatch between the potential energy term involving $x_{k+1}$ and the mixed energy term involving $v_k$ and $x_k$.  Moreover,~\citet{li2024linear} also adopted the mixed sequence  $\{\mathbf{R}_k\}_{k=0}^{\infty}$, but their Lyapunov function includes a kinetic energy term that fundamentally depends on the position update~\eqref{eqn: nag-phase-x}. Therefore, although it achieves linear convergence, it cannot be extended to the \texttt{M-NAG} setting, where the position update is not explicitly available. 

In this study, our key construction is to shift the right-hand side of the previous iterative difference~\eqref{eqn: nag-iv-iter}, $- \left( k+r\right)s \nabla f(y_{k})$, to the left-hand side, and to symmetrically add the new term $- \left( k + r + 1 \right)s\nabla f(y_{k+1})$ to both sides.  This yields a new iterative difference as
\begin{equation}
\left[ \underbrace{\mathbf{R}_{k+1} - \left( k + r + 1 \right)s\nabla f(y_{k+1})}_{:=\mathbf{S}_{k+1}} \right] - \left[ \underbrace{\mathbf{R}_{k} - \left( k+r\right)s \nabla f(y_{k})}_{:=\mathbf{S}_k} \right]  = - \left( k + r + 1 \right)s\nabla f(y_{k+1}) \label{eqn: nag-iv-iter-new},
\end{equation}
where we define a new mixed sequece $\mathbf{S}_k := \mathbf{R}_{k} - \left( k+r\right)s \nabla f(y_{k})$.  A detailed comparision of this new construction with previous formulations of mixed sequences is provided in~\Cref{tab: mix-sequence-iter-diff}. 
\begin{table}[htbp!]
\centering
\begin{tabular}{lll}
\toprule
 & Mixed-Sequence & Iterative Difference \\
\midrule
\parbox[c]{3.5cm}{
\citet{su2016differential}        \\
\citet{shi2022understanding} \\
\citet{chen2022gradient}       \\
\citet{li2022proximal}            \\
\citet{li2024linear} 
}
& \parbox[l]{4.4cm}{$\mathbf{R}_k = (k-1)\sqrt{s}v_{k} + rx_k$ \\ [3pt] ~\eqref{eqn: mix-previous} }
& \parbox[l]{5.8cm}{$\mathbf{R}_{k+1} - \mathbf{R}_k = - (k+r)s \nabla f(y_k)$\\ [3pt] ~\eqref{eqn: nag-iv-iter} } \\
\midrule
\textbf{This Work}
& \parbox[l]{4.4cm}{$\mathbf{S}_k = \mathbf{R}_k  - (k+r)s \nabla f(y_k)$ \\ [3pt]}
& \parbox[l]{5.8cm}{$\mathbf{S}_{k+1} - \mathbf{S}_k = - (k+r+1)s \nabla f(y_{k+1})$\\ [3pt] ~\eqref{eqn: nag-iv-iter-new} }  \\
\bottomrule
\end{tabular}
\caption{Comparison of mixed-sequence constructions and their iterative differences across related literature. Our proposed sequence $\{\mathbf{S}_k\}_{k=0}^{\infty}$  refines the classical formulation by absorbing the gradient term directly, yielding a cleaner iteration with forward-indexed gradients.}
\label{tab: mix-sequence-iter-diff}
\end{table}

As shown in~\Cref{subsubsec: monotonicity}, the mixed-sequece iteration is logically equivalent to the~\texttt{M-NAG} update rule. By substituting the phase-coupling relation~\eqref{eqn: nag-phase-x-y-v} to the definition of $\mathbf{S}_k$, we eliminate the velocity term $v_{k}$, yielding:
\begin{equation}
\label{eqn: mix-new}
\mathbf{S}_k  =  (k-1)\sqrt{s}v_{k} + r x_k - \left( k+r\right)s \nabla f(y_{k}) = (k + r)y_{k}  - kx_k  - (k + r)s\nabla f(y_k).
\end{equation}
Using this expression, the new iterative difference~\eqref{eqn: nag-iv-iter-new} can be reformulated as:
\begin{multline}
(k + r + 1) y_{k+1}  - (k + 1) x_{k+1}  - \left( k + r + 1 \right)\nabla f(y_{k+1}) \\
= (k + r)y_{k}  - kx_k  - (k + r)s\nabla f(y_k) - \left( k + r + 1 \right)\nabla f(y_{k+1}), \label{eqn: m-nag-iteration-new} 
\end{multline}
which is exactly equivalent to the~\texttt{M-NAG} scheme~\eqref{eqn: essential-iteration}. In this study, we also adopt the reformulated iteration~\eqref{eqn: m-nag-iteration-new} as the foundation for constructing the Lyapunov function. Together with the earlier formulation~\eqref{eqn: nag-iv-iter}, the result demonstrate that the construction of a Lyapunov function fundamentally depends only on the \texttt{M-NAG} formulation, without requiring the full update structure of classical~\texttt{NAG}.
\subsection{Overview of contributions}
\label{subsec: contributions}

In this study, building upon the two key observations outlined above, we make the following contributions for the analysis of the forward-backward accelerated algorithms:

\begin{itemize}
\item[(\textbf{I})] In this study, we first highlight that the~\texttt{M-NAG} update is a linear combination of the gradient step and the momentum step in classical~\texttt{NAG}. Since two iterative sequences are embedded in a single update, we gain the flexibility to adjust the auxiliary sequence $\{y_k\}_{k=0}^{\infty}$ to ensure that the objective values $\{f(x_k)\}_{k=0}^{\infty}$ are non-increasing. Moreover, through the implicit-velocity phase-space representation, we find that recovering the full~\texttt{NAG} scheme from~\texttt{M-NAG} still requires one additional assumption, either the position update or the phase-coupling relation, both of which are independently sufficient.

\item[(\textbf{II})] A principled approach to constructing Lyapunov functions is outlined in~\citet{chen2022revisiting}. When the kinetic energy term is excluded, and building on prior works ~\citep{su2016differential, shi2022understanding, chen2022gradient, li2022proximal} for general convex functions, together with our current analysis for strongly convex functions, we demonstrate that the~\texttt{M-NAG} update alone suffices to derive convergence rates. In other words, it is not required to know the full update structure of classical \texttt{NAG}.

\item[(\textbf{III})] In contrast to prior work, we introduce a new mixed sequence $\mathbf{S}_k=\mathbf{R}_k - (k+r)s \nabla f(y_k)$, as defined in~\eqref{eqn: mix-new}, replacing the earlier sequence $\mathbf{R}_k$ from~\eqref{eqn: mix-previous}.  This leads to a new iterative difference, given in~\eqref{eqn: nag-iv-iter-new} or~\eqref{eqn: m-nag-iteration-new}, which facilitates the construction of a novel Lyapunov function. This new Lyapunov function not only guarantees linear convergence for \texttt{NAG}, but by eliminating the kinetic energy term, also enables the generalization of linear convergence to~\texttt{M-NAG}.  Notably, the required starting index admits an explicit expression depending only on the momentum parameter $r$, and is independent of the strong convexity constant, the Lipschitz constant, and the step size. Compared to~\citet{li2024linear}, this yields a more general and robust theoretical result.

\item[(\textbf{IV})]  Finally, leveraging two proximal inequalities for composite functions developed in~\citep{li2024linear, li2024linear2}, which serves as  proximal  proximal analogs of the standard inequalities for strongly convex functions, we further extend our linear convergence to two widely used proximal algorithms: the fast iterative shrinkage-thresholding algorithm~(\texttt{FISTA}) and its monotonic variant~(\texttt{M-FISTA}).

\end{itemize}

\subsection{Related works and organization}
\label{subsec: related}

A major milestone in the history of gradient-based optimization was the introduction of~\texttt{NAG}, a two-step forward-backward algorithm proposed by~\citet{nesterov1983method}. This breakthrough significantly improved convergence rates by incorporating an acceleration mechanism into vanilla gradient descent. Building upon this foundation,~\citet{beck2009fast} introduced a proximal version of the fundamental inequality, leading to the development of~\texttt{FISTA}, a powerful algorithm for composite optimization problems  that has been widely adopted in applications such as signal processing and image reconstruction. The linear convergence of such accelerated algorithms was later established in~\citet{li2024linear}, who employed the high-resolution differential equation framework developed by~\citet{shi2022understanding}. In parallel, ~\citet{beck2009fastieee} proposed a monotonic variant, \texttt{M-FISTA}, which ensures that the objective value is non-increasing across iterations. However, whether linear convergence can be generalized to monotonic variants such as \texttt{M-NAG} and \texttt{M-FISTA} has remained an open question.

Over the past decade, growing interest has focused on understanding acceleration from a dynamical systems perspective. Foundational works~\citep{attouch2012second, attouch2014dynamical} reignited interest in the continuous-time dynamics underlying accelerated methods. This line of inquiry was further advanced by~\citet{su2016differential}, who proposed a low-resolution ODE for modeling \texttt{NAG}, and by~\citep{wibisono2016variational}, who proposed a variational framework to characterize acceleration.  Additional insight came from~\citet{attouch2016rate}, who analyzed the faster convergence rates of objective values. More recently,~\citet{attouch2022ravine} investigated new mixed sequences to unify the treatment of the Ravine method and~\texttt{NAG}. However, their study did not address monotonic variants, nor did it explore how such variants, such as \texttt{M-NAG}, are related to classical \texttt{NAG} in terms of structure or convergence behavior. A complete understanding of the acceleration mechanism emerged through comparisons between \texttt{NAG} and Polyak’s heavy-ball method, with pivotal contributions from the high-resolution differential equation framework proposed by~\citet{shi2022understanding}. This framework revealed that acceleration arises from gradient correction. Follow-up works~\citep{chen2022gradient, chen2022revisiting} demonstrated that this implicit-velocity perspective is particularly effective for analyzing \texttt{NAG}. Subsequent advances extended this framework to composite optimization, covering both convex and strongly convex settings~\citep{li2022proximal, li2024linear2},  and further to overdamped regimes~\citep{chen2023underdamped}. These developments have substantially deepened our understanding of acceleration and paved the way for designing more efficient and theoretically grounded optimization algorithms.

The remainder of this paper is organized as follows.~\Cref{sec: prelim} introduces the foundational definitions and inequalities for strongly convex functions and composite functions, which serve as preliminaries for the analysis.~\Cref{sec: nag} outline the construction of the Lyapunov function and establishes the linear convergence for~\texttt{NAG} and~\texttt{FISTA}.~\Cref{sec: m-nag-fista} extends the linear convergence to the monotonically accelerated forward-backward algorithms,~\texttt{M-NAG} and \texttt{M-FISTA}, leveraging the novel Lyapnov function that excludes the kinetic energy. Finally,~\Cref{sec: conclusion} concludes this papers and proposes potential avenues for future research.

\section{Preliminaries}
\label{sec: prelim} 

In this paper, we adopt the notations that are consistent with those used in~\citet{nesterov2018lectures},~\citet{shi2022understanding}, and~\citet{li2024linear2}, with slight adjustments tailored to the context of our analysis. Let $\mathcal{F}^0(\mathbb{R}^d)$ denote the class of continuous convex functions defined on $\mathbb{R}^d$. Specifically, $g \in \mathcal{F}^{0}(\mathbb{R}^d)$ if it satisfies the convex condition:
\[
g\left( \alpha x + (1 - \alpha)y \right) \leq \alpha g(x) + (1 - \alpha)g(y),
\] 
for any $x,y \in \mathbb{R}^d$ and $\alpha \in [0, 1]$. A subclass $\mathcal{F}^0(\mathbb{R}^d) \subset \mathcal{F}^1_L(\mathbb{R}^d)$ consists of continuously differentiable functions with Lipschitz continuous gradients. Specifically, $f \in \mathcal{F}^{1}_{L}(\mathbb{R}^d)$ if $f \in \mathcal{F}^{0}(\mathbb{R}^d)$ and its gradient satisfies:
\begin{equation}
\label{eqn: grad-lip}
\| \nabla f(x) - \nabla f(y) \| \leq L \| x - y \|,
\end{equation}
for any $x, y \in \mathbb{R}^d$. Next, we define the class $\mathcal{S}_{\mu,L}^{1}(\mathbb{R}^d) \subset \mathcal{F}^1_L(\mathbb{R}^d)$, consisting of functions that are $\mu$-strongly convex for some $0 < \mu \leq L$. Specifically, $f \in \mathcal{S}_{\mu,L}^{1}(\mathbb{R}^d)$ if $f \in \mathcal{F}^1_L(\mathbb{R}^d)$ and satisfies the $\mu$-strongly convexity condition:
\begin{equation}
\label{eqn: defn-scvx}
f(y) \geq f(x) +  \langle \nabla f(y), y - x \rangle + \frac{\mu}{2} \|y - x\|^2,
\end{equation}
for any $x, y \in \mathbb{R}^d$. For any $f \in \mathcal{S}_{\mu,L}^{1}(\mathbb{R}^d)$, the following fundamental inequality holds:
\begin{equation}
\label{eqn: fund-inq-smooth}
f(y - s\nabla f(y)) - f(x) \leq \left\langle \nabla f(y), y - x \right\rangle - \frac{\mu}{2} \|y - x\|^2- \left( s - \frac{Ls^2}{2} \right) \|\nabla f(y)\|^2,
\end{equation}
for any $x,y \in \mathbb{R}^d$. Additionally, let $x^\star$ as the unique minimizer of the objective function $f$. Then, for any $f \in \mathcal{S}_{\mu,L}^{1}(\mathbb{R}^d)$ and any $y \in \mathbb{R}^d$, the following inequality holds: 
\begin{equation}
\label{eqn: key-inq-smooth}
\left\| \nabla f(y) \right\|^2 \geq 2\mu \left( f(y - s\nabla f(y)) - f(x^{\star}) \right).
\end{equation}

We now turn to the composite objective $\Phi= f + g$, where $f \in \mathcal{S}_{\mu,L}^1(\mathbb{R}^d)$ and $g \in \mathcal{F}^0(\mathbb{R}^d)$. Following~\citet{beck2009fast} and~\citet{su2016differential}, we introduce the notion of the $s$-proximal value, a central concept in our analysis.
\begin{defn}[$s$-proximal value]
\label{defn: proximal-value}
Let the step size satisfy $s \in (0, 1/L)$. For any $f\in \mathcal{S}_{\mu,L}^1(\mathbb{R}^d)$ and $g \in \mathcal{F}^0(\mathbb{R}^d)$, the $s$-proximal value is defined as
\begin{equation}
\label{eqn: proximal-value}
    P_s(x) := \mathop{\arg\min}_{y\in\mr^d}\left\{ \frac{1}{2s}\left\| y - \left(x - s\nabla f(x)\right) \right\|^2 + g(y) \right\}.
\end{equation}
The $s$-proximal value minimizes a weighted sum of the squared distance from the gradient step $x - s\nabla f(x)$, plus the regularization term $g(y)$. The formulation is especially powerful for handling composite objectives with non-smooth but proximally tractable regularizers. For instance, when $g(x) = \lambda \|x\|_1$, the $s$-proximal value admits a closed-form expression. Specifically, for any $x \in \mathbb{R}^d$, the $i$-th component of $s$-proximal value is given by: 
\[
P_s(x)_i = \big(\left|\left(x - s\nabla f(x)\right)_i\right| - \lambda s\big)_+ \text{sgn}\big(\left(x - s\nabla f(x)\right)_i\big), 
\]
for $i=1,\ldots,d$, where $\text{sgn}(\cdot)$ denotes the sign function. 
\end{defn}

\begin{defn}[$s$-proximal subgradient]
\label{defn: subgradient}
The $s$-proximal subgradient at  $x \in \mathbb{R}^d$ is defined as:
\begin{equation}
\label{eqn: subgradient}
G_s(x): = \frac{x - P_s(x)}{s}.
\end{equation}
This quantity generalizes the notion of a gradient to the composite setting, capturing the direction of descent for $\Phi = f + g$ when the proximal step is applied.

\end{defn}

Using this definition,  we extend the two key inequalities,~\eqref{eqn: fund-inq-smooth} and~\eqref{eqn: key-inq-smooth}, to the proximal setting.

\begin{lemma}[Lemma 4 in~\citet{li2024linear2}]
\label{lem: fund-inq-composite}
Let $\Phi = f + g$ be a composite function with $f \in \mathcal{S}_{\mu, L}^1(\mathbb{R}^d)$ and $g \in \mathcal{F}^0(\mathbb{R}^d)$. Then, for any $x, y \in \mathbb{R}^d$, the following inequality holds:
\begin{equation}
\label{eqn: fund-inq-composite}
\Phi(y - sG_s(y)) - \Phi(x) \leq \left\langle G_s(y), y - x \right\rangle - \frac{\mu}{2}\|y - x\|^2 - \left( s - \frac{Ls^2}{2} \right) \|G_s(y)\|^2.
\end{equation}

\end{lemma}
\begin{lemma}[Lemma 4.3 in~\citet{li2024linear}]
\label{lem: key-inq-composite}
Let $f \in \mathcal{S}_{\mu,L}^{1}(\mathbb{R}^d)$ and $g \in \mathcal{F}^0(\mathbb{R}^d)$. Then, for any $ y \in \mathbb{R}^d$, the s-proximal subgradient satisfies the following inequality:
\begin{equation}
\label{eqn: key-inq-composite}
\|G_s(y)\|^2 \geq 2\mu \left( \Phi(y - sG_s(y)) - \Phi(x^{\star}) \right),
\end{equation}
where we shall note that here $x^{\star}$ is the unique minimizer of the objective function $\Phi$.
\end{lemma}

\section{Lyapnov analysis for NAG and FISTA}
\label{sec: nag}

In this section, we begin by constructing a novel Lyapunov function, following the principled approach proposed in~\citet{chen2022revisiting}. Central to this construction is the mixed sequence $\mathbf{S}_{k}$, defined in~\eqref{eqn: mix-new}, which forms the foundation of our Lyapunov analysis. Building on this, we establish the linear convergence rate of~\texttt{NAG} for strongly convex functions and present the result in a formal theorem. We then apply~\Cref{lem: fund-inq-composite} and~\Cref{lem: key-inq-composite} to extend the linear convergence guarantee to the proximal setting, with particular emphasis on the~\texttt{FISTA} algorithm. 

\subsection{The construction of a novel Lyapunov function}
\label{subsec: construction}

Following the principled approach outlined in~\citet{chen2022revisiting}, we construct a novel Lyapunov function in two main steps: first, by introducing the new mixed sequence $\mathbf{S}_{k}$ to define the mixed energy, and second, by adjusting the coefficient of the potential energy to align with this mixed energy.

\begin{itemize}
\item[(\textbf{1})]

\textbf{Construction of the mixed energy} Recall the new mixed sequence $\mathbf{S}_k$ defined in~\eqref{eqn: mix-new}. In convergence analysis, the quantity of interest is not the iterate $x_k$ itself, but its difference from the optimal solution $x_k - x^{\star}$. To incorporate this, we introduce the term $-rx^{\star}$ to $\mathbf{S}_k$ and $\mathbf{S}_{k+1}$ in the iterative difference~\eqref{eqn: nag-iv-iter-new}, resulting in:
\begin{multline}
k\sqrt{s}v_{k+1} + r ( x_{k+1} - x^{\star} ) - \left( k + r + 1 \right)s\nabla f(y_{k+1})  \\
= (k-1)\sqrt{s}v_{k} + r ( x_k - x^{\star} ) - \left( k+r\right)s \nabla f(y_{k}) - \left( k + r + 1 \right)s\nabla f(y_{k+1}).   \label{eqn: nag-phase-iteration-new-xstar}
\end{multline}
This relation naturally leads to the mixed energy, which is given by: 
\begin{equation}
\label{eqn: mix}
\mathcal{E}_{mix}(k) = \frac12 \left\| (k-1)\sqrt{s}v_{k} + r (x_{k} - x^{\star}) - s (k + r) \nabla f(y_k) \right\|^2.
\end{equation}
Next, using the iterative relation~\eqref{eqn: nag-phase-iteration-new-xstar}, we analyze the difference between successive mixed energy iterates:
\begin{align}
\mathcal{E}&_{mix}(k + 1) - \mathcal{E}_{mix}(k) \nonumber \\ 
                  & = \left\langle - \left( k + r + 1 \right)s\nabla f(y_{k+1}), k\sqrt{s}v_{k+1} + r ( x_{k+1} - x^{\star} ) - \frac{\left( k + r + 1 \right)s}{2} \cdot \nabla f(y_{k+1}) \right\rangle. \label{eqn: mix-iter1}
\end{align}
Substituting~\eqref{eqn: nag-phase-x-y-v} into~\eqref{eqn: mix-iter1}, we derive the following result for the iterative difference:  
\begin{align}
\mathcal{E}_{mix}(k + 1) - \mathcal{E}_{mix}(k)  = & \underbrace{- k(k+1)s\sqrt{s} \left\langle \nabla f(y_{k+1}), v_{k+1}\right\rangle}_{\mathbf{I}} \nonumber \\
                                                                                & - r (k + r + 1) s \left\langle \nabla f(y_{k+1}), y_{k + 1} - x^{\star} \right\rangle \nonumber \\
                                                                                & + \frac{(k+r+1)^2s^2}{2}\left\| \nabla f(y_{k+1}) \right\|^2. \label{eqn: mix-iter2}
\end{align}

\item[(\textbf{2})] \textbf{Construction of potential energy} In contrast to the Lyapunov function defined in~\eqref{eqn: lyapunov-cvx}, previously used for convex functions~\citep{chen2022gradient, li2022proximal}, we now modify the potential energy by replacing $f(x_{k})$ with $f(x_{k+1})$ to better align with the mixed energy in~\eqref{eqn: mix}. Additionally, we introduce a dynamic coefficient $\tau(k)$ to the potential energy, defined as: 
\begin{equation}
\label{eqn: pot}
\mathcal{E}_{pot}(k) := s \tau(k) \left( f(x_{k+1}) - f(x^{\star}) \right).
\end{equation}
The iterative difference of this potential energy is then given by:
\begin{align}
\mathcal{E}_{pot}(k+1) &- \mathcal{E}_{pot}(k) \nonumber \\ & = s \tau(k) \left( f(x_{k+2}) - f(x_{k+1}) \right) + s\left( \tau(k+1) - \tau(k) \right) \left( f(x_{k+2}) - f(x^{\star}) \right).    \label{eqn: pot-iter1}
\end{align}
By applying the fundamental inequality~\eqref{eqn: fund-inq-smooth} into $x_{k+2}$ and $x_{k+1}$, we have
\begin{align}
 f(x_{k+2}) - f(x_{k+1})  \leq &  \left\langle \nabla f(y_{k+1}), y_{k+1} - x_{k+1}\right\rangle \nonumber \\
                                            & - \frac{\mu}{2} \left\| y_{k+1} - x_{k+1} \right\|^2 - \left( s - \frac{Ls^2}{2} \right) \left\| \nabla f(y_{k+1}) \right\|^2  \label{eqn: fund-inq-xkyk-nag}
\end{align}
Substituting~\eqref{eqn: fund-inq-xkyk-nag} into~\eqref{eqn: pot-iter1}, we obtain the following bound on the potential energy difference:
\begin{align}
\mathcal{E}&_{pot}(k+1) - \mathcal{E}_{pot}(k) \nonumber \\
                  & \leq s \tau(k) \left( \left\langle \nabla f(y_{k+1}), y_{k+1} - x_{k+1} \right\rangle - \frac{\mu}{2} \left\| y_{k+1} - x_{k+1} \right\|^2 - \left( s - \frac{Ls^2}{2} \right) \left\| \nabla f(y_{k+1}) \right\|^2 \right) \nonumber \\
                  & \mathrel{\phantom{\leq}} + s \left( \tau(k+1) - \tau(k) \right) \left( f(x_{k+2}) - f(x^{\star}) \right).
                   \label{eqn: pot-iter2} 
\end{align}      
Substituting~\eqref{eqn: nag-phase-x} into~\eqref{eqn: pot-iter2}, we further refine the inequality as:            
\begin{align}    
\mathcal{E}&_{pot}(k+1) - \mathcal{E}_{pot}(k) \nonumber \\              
                  & =    s \tau(k)  \left(  \underbrace{\frac{k\sqrt{s}\left\langle \nabla f(y_{k+1}), v_{k+1}  \right\rangle }{k+r+1} }_{\mathbf{II}} - \frac{\mu}{2} \cdot \frac{k^2 s\left\| v_{k+1} \right\|^2}{(k+r+1)^2} - \left( s - \frac{Ls^2}{2} \right) \left\| \nabla f(y_{k+1}) \right\|^2 \right)  \nonumber \\
                  & \mathrel{\phantom{\leq}} + s \left( \tau(k+1) - \tau(k) \right) \left( f(x_{k+2}) - f(x^{\star}) \right)  \label{eqn: pot-iter3} 
\end{align}
To  ensure cancellation of terms $\mathbf{I} + (k+1)(k+r+1)s \cdot \mathbf{II} = 0$, we choose the dynamic coefficient as  $\tau(k) = (k+1)(k+r+1)$. 
\end{itemize}

Combining the mixed energy from~\eqref{eqn: mix} and the potential energy from~\eqref{eqn: pot}, the Lyapunov function is constructed as: 
\begin{align}
\mathcal{E}(k) = &\; s(k+1)(k+r+1)\left( f(x_{k+1}) - f(x^{\star}) \right) \nonumber \\ 
                           &\; + \frac12\left\| (k-1)\sqrt{s}v_{k} + r (x_{k} - x^{\star}) - s (k + r) \nabla f(y_k) \right\|^2. \label{eqn: lyapunov-nag}
\end{align}
This Lyapunov function encapsulates both the mixed and potential energies and plays a key role in analyzing the linear convergence of the iterative process.

\subsection{Linear convergence of NAG}
\label{subsec: lin-nag}
Using the Lyapunov function defined in~\eqref{eqn: lyapunov-nag}, we now establish the linear convergence of~\texttt{NAG}. The result is stated formally in the following theorem:

\begin{theorem}
\label{thm: nag}
Let $f \in \mathcal{S}_{\mu,L}^{1}(\mathbb{R}^d)$. Given any step size $0 < s < 1/L$, there exists a positive integer $K: = \max\left\{0,  \frac{3r^2 - 4r - 12}{8}\right\}$ such that the iterative sequence $\{x_{k}\}_{k=0}^{\infty}$ generated by~\texttt{NAG}~\eqref{eqn: nag}, with any initial point $x_0 = y_0 \in \mathbb{R}^d$, satisfies the following inequality
\begin{equation}
\label{eqn: nag-rate}
f(x_k) - f(x^{\star}) \leq \frac{(r +1)\left( f(x_1) - f(x^{\star})\right) + r^2 L\| x_1 - x^{\star} \|^2}{k(k+r) \left[ 1 + (1 - Ls) \cdot \frac{\mu s}{4} \right]^k}, 
\end{equation}
for any $k \geq \max\left\{1,K \right\}$. 
\end{theorem}

\begin{proof}[Proof of~\Cref{thm: nag}]
We begin by deriving a decrement of the Lyapunov function $\mathcal{E}(k)$ across successive iterations. By combining the expressions for the mixed and potential energy differences in~\eqref{eqn: mix-iter2} and~\eqref{eqn: pot-iter3}, we obtain the following bound for the Lyapunov difference:
\begin{align}
\mathcal{E}(k+1) - \mathcal{E}(k) = & - r (k + r + 1) s \left\langle \nabla f(y_{k+1}), y_{k + 1} - x^{\star} \right\rangle  + \frac{(k+r+1)^2s^2}{2}\left\| \nabla f(y_{k+1}) \right\|^2 \nonumber \\ 
                                                         & - \frac{\mu}{2} \cdot \frac{s^2k^2(k+1)}{k+r+1} \cdot \|v_{k+1}\|^2 - s(k+1)(k+r+1)\left( s - \frac{Ls^2}{2} \right) \| \nabla f(y_{k+1}) \|^2  \nonumber \\
                                                         & + s (2k + r + 3) \left( f(x_{k+2}) - f(x^{\star}) \right). \label{eqn: lyapunov-iter1}
\end{align}
To estimate the first inner product term,  we apply the fundamental inequality~\eqref{eqn: fund-inq-smooth} with $x_{k+2}$ and $x^{\star}$, yielding:  
\begin{align}
- \left\langle \nabla f(y_{k+1}), y_{k+1} - x^{\star}\right\rangle \leq & - \left( f(x_{k+2}) - f(x^{\star}) \right) \nonumber \\ &- \frac{\mu}{2} \left\| y_{k+1} - x^{\star} \right\|^2 - \left( s - \frac{Ls^2}{2} \right) \left\| \nabla f(y_{k+1}) \right\|^2. \label{eqn: fund-inq-xstar}
\end{align}
Substituting~\eqref{eqn: fund-inq-xstar} into~\eqref{eqn: lyapunov-iter1}, we obtain:
\begin{align}
\mathcal{E}(k+1) - \mathcal{E}(k) \leq & - s \left[(r -2)k + (r^2 - 3) \right] \left( f(x_{k+2}) - f(x^{\star}) \right)  - \frac{\mu}{2} \cdot \frac{s^2k^2(k+1)}{k+r+1} \cdot \|v_{k+1}\|^2     \nonumber \\
                                                             & - \frac{ s (2k + r + 3) \mu}{2} \|y_{k+1} - x^{\star}\|^2  - \frac{(k+r+1)^2s^2(1-Ls)}{2}\left\| \nabla f(y_{k+1}) \right\|^2.           \label{eqn: lyapunov-iter2} 
\end{align}
To simplify further, we apply the strong convexity inequality~\eqref{eqn: key-inq-smooth}, which gives a lower bound on the gradient norm as: 
\begin{equation}
\left\| \nabla f(y_{k+1}) \right\|^2 \geq 2\mu \left( f(x_{k+2}) - f(x^{\star}) \right). 
\label{eqn: fund-inq-smooth1-new}
\end{equation}
Substituting~\eqref{eqn: fund-inq-smooth1-new} into~\eqref{eqn: lyapunov-iter2},  we arrive at:
\begin{align}
\mathcal{E}(k+1) - \mathcal{E}(k) \leq & - \frac{\mu s(1-Ls)}{4} \cdot s(k+r+1)^2\left( f(x_{k+2}) - f(x^{\star}) \right) \nonumber \\ & - \frac{\mu s^2}{2} \cdot \frac{k^2(k+1)}{k+r+1} \cdot \|v_{k+1}\|^2     \nonumber \\
                                                             & - \frac{ s (2k + r + 3) \mu}{2} \|y_{k+1} - x^{\star}\|^2  - \frac{3(k+r+1)^2s^2(1-Ls)}{8}\left\| \nabla f(y_{k+1}) \right\|^2.           \label{eqn: lyapunov-iter3} 
\end{align}
Next, we estimate the Lyapunov function itself using the Cauchy-Schwarz inequality. From the definition~\eqref{eqn: lyapunov-nag}, we have:
\begin{align}
\mathcal{E}(k + 1) \leq  & s(k+2)(k+r+2)\left( f(x_{k+2}) - f(x^{\star}) \right) + \frac{3s}{2} \cdot \frac{k^2(k+1)^2}{(k+r+1)^2} \cdot \left\| v_{k+1} \right\|^2 \nonumber  \\
                                      & + \frac{3r^2}{2} \left\| y_{k+1} - x^{\star} \right\|^2 + \frac{3s^2}{2} \cdot (k+r+1)^2 \cdot \|\nabla f(y_{k+1}) \|^2.                         \label{eqn: lyapunov-estimate} 
\end{align}
By comparing the coefficients of the corresponding terms in inequalities~\eqref{eqn: lyapunov-iter3} and~\eqref{eqn: lyapunov-estimate}, we establish a linear contraction of the Lyapunov function. Specifically, we establish the bound
\begin{align}
\mathcal{E}(k+1) - \mathcal{E}(k) & \leq - \mu s \cdot \min\left\{ \frac{1-Ls}{4}, \frac{1}{3}, \frac{2k+r+3}{3r^2}, \frac{1-Ls}{4 \mu s} \right\} \mathcal{E}(k+1)    \nonumber \\
                                                     & \leq - \mu s \cdot \left(\frac{1 - Ls}{4} \right) \cdot \mathcal{E}(k+1),                                                                                     \label{eqn: iter-diff-lyapunov} 
\end{align}
where the last inequality follows from the choice of $K$ and the fact that $0 < \mu s < \mu/L \leq 1$. This concludes the proof by invoking the Lipschitz gradient condition~\eqref{eqn: grad-lip}. 
\end{proof}


\subsection{Linear convergence of FISTA}
\label{subsec: lin-fista}

We now extend the linear convergence established in~\Cref{thm: nag} for~\texttt{NAG} to its proximal variant,~\texttt{FISTA}. As specified in~\Cref{defn: proximal-value},~\texttt{FISTA} utilizes the $s$-proximal value~\eqref{eqn: proximal-value} and proceeds via the following iterative scheme, starting from any initial point $y_0 = x_0 \in \mathbb{R}^d$:
\begin{subequations}
\label{eqn: fista}
\begin{empheq}[left=\empheqlbrace]{align}
& x_{k+1} = P_s\left( y_{k} \right),                                                                           \label{eqn: fista-gradient}          \\
& y_{k+1} = x_{k+1} + \frac{k}{k+r+1} (x_{k+1} - x_{k}),                                          \label{eqn: fista-momentum} 
\end{empheq}    
\end{subequations}
where $s>0$ is the step size. According to~\Cref{defn: subgradient}, this update rule can be equivalently written in terms of $s$-proximal subgradient $G_s(y_k)$ (see~\eqref{eqn: subgradient}), yielding a form analogous to~\texttt{NAG}: 
\begin{subequations}
\label{eqn: fista_1}
\begin{empheq}[left=\empheqlbrace]{align}
& x_{k+1} = y_{k} - sG_s(y_{k}),                                                                              \label{eqn: fista_1-gradient}          \\
& y_{k+1} = x_{k+1} + \frac{k}{k+r+1} (x_{k+1} - x_{k}),                                          \label{eqn: fista_1-momentum} 
\end{empheq}    
\end{subequations}
where $G_s(y_{k})$ serves as a natural generalization of the classical gradient $\nabla f(y_k)$. This formulation highlights the structural similarity between~\texttt{NAG} and~\texttt{FISTA}, while accounting for the non-smooth nature of $\Phi$. To further parallel the \texttt{NAG} formulation, we define the velocity sequence as$v_{k}:= (x_{k}- x_{k-1})/\sqrt{s}$. This allows us to rewrite the \texttt{FISTA} iterations in the implicit-velocity phase-space form:
\begin{subequations}
\label{eqn: fista-phase}
\begin{empheq}[left=\empheqlbrace]{align}
         & x_{k+1} - x_{k} = \sqrt{s} v_{k+1},                                                                               \label{eqn: fista-phase-x}          \\
         & v_{k+1} - v_{k} = -\frac{r+1}{k+r}  v_{k} - \sqrt{s} G_s \left( y_{k} \right),                 \label{eqn: fista-phase-v}
\end{empheq}    
\end{subequations}
where the sequence $\{y_{k}\}_{k=0}^{\infty}$ satisfies the phase-coupling relationship:
\begin{equation}
\label{eqn: fista-phase-x-y-v}
y_k = x_{k} +  \frac{k-1}{k+r} \cdot \sqrt{s}v_k.
\end{equation}

To prove linear convergence of~\texttt{FISTA}, we construct a generalized Lyapunov function by modifying the one in~\eqref{eqn: lyapunov-nag}, replacing the smooth objective $f$ with the composite objective $\Phi = f+g$. This yields the following Lyapunov function: 
\begin{align}
\mathcal{E}(k) = &\; s(k+1)(k+r+1)\left( \Phi(x_{k+1}) - \Phi(x^{\star}) \right) \nonumber \\ 
                           &\; + \frac12\left\| (k-1)\sqrt{s}v_{k} + r (x_{k} - x^{\star}) - s (k + r) G_s(y_k) \right\|^2. \label{eqn: lyapunov-fista}
\end{align}
As seed in~\Cref{subsec: lin-nag}, the linear convergence of~\texttt{NAG} rests on three key inequalities,~\eqref{eqn: fund-inq-xkyk-nag},~\eqref{eqn: fund-inq-xstar}, and~\eqref{eqn: fund-inq-smooth1-new}, which characterize how the iterates approach optimality. In the proximal setting, these inequalities are generalized via~\Cref{lem: fund-inq-composite} and~\Cref{lem: key-inq-composite}, which extend the smooth-case strong convexity results~\eqref{eqn: fund-inq-smooth} and~\eqref{eqn: key-inq-smooth} to accommodate the non-smooth term $g$. These extensions allow us to adapt the analysis of~\texttt{NAG} directly to~\texttt{FISTA}. We now state the linear convergence result for~\texttt{FISTA}:

\begin{theorem}
\label{thm: fista}
Let $\Phi = f + g$, where $f \in \mathcal{S}_{\mu,L}^{1}(\mathbb{R}^d)$ and $g \in \mathcal{F}^0(\mathbb{R}^d)$. Given any step size $0 < s < 1/L$, there exists a positive integer $K: = \max\left\{0,  \frac{3r^2 - 4r - 12}{8}\right\}$ such that the iterative sequence $\{x_{k}\}_{k=0}^{\infty}$ generated by~\texttt{FISTA}~\eqref{eqn: fista},  with any initial $x_0 = y_0 \in \mathbb{R}^d$, satisfies the following inequality: 
\begin{equation}
\label{eqn: fista-rate}
\Phi(x_k) - \Phi(x^{\star}) \leq \frac{(r +1)\left( \Phi(x_1) - \Phi(x^{\star})\right) + r^2 L\| x_1 - x^{\star} \|^2}{k(k+r) \left[ 1 + (1 - Ls) \cdot \frac{\mu s}{4} \right]^k}, 
\end{equation}
for any $k \geq \max\left\{ 1, K \right\}$. 
\end{theorem}

\section{Monotonicity and linear convergence}
\label{sec: m-nag-fista}

In this section, we revisit the linear convergence results using the new mixed sequence $\mathbf{S}_k$ defined in~\eqref{eqn: mix-new}, this time without relying on the phase-space representation. Rather than employing the velocity iterates $v_k$, we work directly with the auxiliary iterates $y_{k}$. By leveraging the iterative difference~\eqref{eqn: m-nag-iteration-new}, we extend the linear convergence for~\texttt{NAG} established in~\Cref{sec: nag} to its monotonic variant,~\texttt{M-NAG} (see~\eqref{eqn: m-nag} or~\eqref{eqn: essential-iteration}). This demonstrates that the linear convergence rate, as captured by the Lyapunov function, fundamentally relies on the \texttt{M-NAG} update rule~\eqref{eqn: essential-iteration}, rather than the full information from the original \texttt{NAG} iteration. We then further extend this analysis of the linear convergence to cover the proximal monotonic variant, \texttt{M-FISTA}.

\subsection{Smooth optimization via M-NAG}
\label{subsec: m-nag}

Using the mixed sequence $\mathbf{S}_k$ defined in~\eqref{eqn: mix-new}, we write down a new form of the Lyapunov function equivalent to~\eqref{eqn: lyapunov-nag}, which depends only on the iterates, $x_k$ and $y_{k}$: 
\begin{align}
\mathcal{E}(k) = &\;  \underbrace{s(k+1)(k+r+1)\left( f(x_{k+1}) - f(x^{\star}) \right)}_{:=\mathcal{E}_{pot}(k)} \nonumber \\ 
                           &\; +  \underbrace{\frac12\left\| k(y_{k} -x_k) + r(y_k - x^{\star}) - (k + r)s\nabla f(y_k) \right\|^2}_{:=\mathcal{E}_{mix}(k)}.   \label{eqn: lyapunov-m-nag}
\end{align}
This Lyapunov function p is the key to proving the linear convergence of~\texttt{M-NAG}, as stated in the following theorem:
\begin{theorem}
\label{thm: m-nag}
Let $f \in \mathcal{S}_{\mu,L}^{1}(\mathbb{R}^d)$. Given any step size $0 < s < 1/L$, there exists a positive integer $K: = \max\left\{0,  \frac{3r^2 - 4r - 12}{8}\right\}$ such that the iterative sequence $\{x_{k}\}_{k=0}^{\infty}$ generated by~\texttt{M-NAG}~\eqref{eqn: m-nag}, with any initial point $x_0 = y_0 \in \mathbb{R}^d$, satisfies the following inequality
\begin{equation}
\label{eqn: mnag-rate}
f(x_k) - f(x^{\star}) \leq \frac{(r +1)\left( f(x_1) - f(x^{\star})\right) + r^2 L\| x_1 - x^{\star} \|^2}{k(k+r) \left[ 1 + (1 - Ls) \cdot \frac{\mu s}{4} \right]^k}, 
\end{equation}
for any $k \geq \max\left\{1,K \right\}$. 
\end{theorem}

\begin{proof}[Proof of~\Cref{thm: m-nag}]
We compute the iterative difference of the Lyapunov function $\mathcal{E}(k)$, separating it into two components: the potential energy and the mixed energy. 
\begin{itemize}
\item[(\textbf{1})] Using the definition of the potential energy $\mathcal{E}_{pot}(k)$, we compute the following difference:
\begin{align}
 \mathcal{E}_{pot}(k+1) - \mathcal{E}_{pot}(k)  =  & s (k+1)(k+r+1) \left( f(x_{k+2}) - f(x_{k+1}) \right) \nonumber \\
                                                                              & + s \left( 2k + r + 3 \right) \left( f(x_{k+2}) - f(x^{\star}) \right).              \label{eqn: pot-iter-mnag-1}
\end{align}
Next, we apply the fundamental inequality~\eqref{eqn: fund-inq-smooth} to the points $x_{k+2}$ and $x_{k+1}$, and note from~\eqref{eqn: m-nag-monetone} that $ f(x_{k+2}) \leq f(z_{k+1}) = f(y_{k+1} - s \nabla f(y_{k+1})) $. This yields:
\begin{align}
 f(x_{k+2}) - f(x_{k+1})  \leq &  \left\langle \nabla f(y_{k+1}), y_{k+1} - x_{k+1}\right\rangle \nonumber \\
                                            & - \frac{\mu}{2} \left\| y_{k+1} - x_{k+1} \right\|^2 - \left( s - \frac{Ls^2}{2} \right) \left\| \nabla f(y_{k+1}) \right\|^2.  \label{eqn: fund-inq-xkyk-mnag}
\end{align}
Substituting~\eqref{eqn: fund-inq-xkyk-mnag} into~\eqref{eqn: pot-iter-mnag-1}, we obtain the following upper bound:
\begin{align}
\mathcal{E}_{pot}(k+1) - \mathcal{E}_{pot}(k)  \leq & \underbrace{s (k+1)(k+r+1)  \left\langle \nabla f(y_{k+1}), y_{k+1} - x_{k+1} \right\rangle}_{\mathbf{I}}      \nonumber \\
                                                                                 & - \frac{\mu s}{2} \cdot (k+1)(k+r+1) \| y_{k+1} - x_{k+1} \|^2                                                                        \nonumber \\
                                                                                 & - s  \left( s - \frac{Ls^2}{2} \right) (k+1)(k+r+1)  \| \nabla f(y_{k+1}) \|^2                                                       \nonumber \\
                                                                                 & + s \left( 2k + r + 3\right) \left( f(x_{k+2}) - f(x^{\star}) \right).                                                                        \label{eqn: pot-iter-mnag-2}
\end{align}

%
%

\item[(\textbf{2})]
To relate the iterates to the optimal point $x^{\star}$, we introduce the term $-rx^{\star}$ to both sides of the new mixed-sequence iteration~\eqref{eqn: m-nag-iteration-new}, resulting in:
\begin{multline}
(k+1)(y_{k+1} - x_{k+1}) + r ( y_{k+1} - x^{\star} ) - \left( k + r + 1 \right)s\nabla f(y_{k+1})  \\
= k(y_{k} - x_{k}) + r ( y_k - x^{\star} ) - \left( k+r\right)s \nabla f(y_{k}) - \left( k + r + 1 \right)s\nabla f(y_{k+1}).   \label{eqn: nag-iteration-new-xstar}
\end{multline}
Using this identity, we derive the mixed energy difference: 
\begin{align}
\mathcal{E}_{mix}(k+1) - \mathcal{E}_{mix}(k)  = & -s(k+1)(k+r+1) \left\langle \nabla f(y_{k+1}), y_{k+1} - x_{k+1} \right\rangle \nonumber   \\
                                                                              & - sr(k+r+1) \left\langle \nabla f(y_{k+1}), y_{k+1} - x^{\star} \right\rangle                                                   \nonumber   \\
                                                                              &  + \frac{s^2(k+r+1)^2}{2} \| \nabla f(y_{k+1}) \|^2.                                                                                      \label{eqn: mix-iter-mnag-1} 
\end{align}
Now, applying the fundamental inequality~\eqref{eqn: fund-inq-smooth} to $x_{k+2}$ and $x^{\star}$, and rearranging, we obtain: 
\begin{align}
- \left\langle \nabla f(y_{k+1}), y_{k+1} - x^{\star}\right\rangle \leq & - \left( f(x_{k+2}) - f(x^{\star}) \right) \nonumber \\
                                                                                                       & - \frac{\mu}{2} \left\| y_{k+1} - x^{\star} \right\|^2 - \left( s - \frac{Ls^2}{2} \right) \left\| \nabla f(y_{k+1}) \right\|^2.  \label{eqn: fund-inq-mnag}
\end{align}
Substituting~\eqref{eqn: fund-inq-mnag} into~\eqref{eqn: mix-iter-mnag-1}, we get:
\begin{align}
\mathcal{E}_{mix}(k+1) - \mathcal{E}_{mix}(k)  \leq & - \underbrace{s(k+1)(k+r+1) \left\langle \nabla f(y_{k+1}), y_{k+1} - x_{k+1} \right\rangle}_{\mathbf{II}}   \nonumber   \\
                                                                                  & - sr(k+r+1)  \left( f(x_{k+2}) - f(x^{\star}) \right)                                                                                           \nonumber   \\
                                                                                  & -  \frac{\mu s}{2} \cdot  r(k+r+1)\left\| y_{k+1} - x^{\star} \right\|^2                                                              \nonumber   \\
                                                                                  &  - s\left( s - \frac{Ls^2}{2} \right) \cdot r(k+r+1) \left\| \nabla f(y_{k+1}) \right\|^2                                        \nonumber   \\
                                                                                  &  + \frac{s^2(k+r+1)^2}{2} \| \nabla f(y_{k+1}) \|^2.                                                                                      \label{eqn: mix-iter-mnag-2} 
\end{align}
\end{itemize}

From the iterative differences of the potential and the mixed energies in~\eqref{eqn: pot-iter-mnag-2} and~\eqref{eqn: mix-iter-mnag-2}, we observe that $\mathbf{I} = \mathbf{II}$. This allows us to combine the two expressions and obtain the following bound on the Lyapunov difference:
\begin{align}
\mathcal{E}(k+1) - \mathcal{E}(k)  \leq  & - s \left[(r -2)k + (r^2 - 3) \right] \left( f(x_{k+2}) - f(x^{\star}) \right)                     \nonumber                    \\
                                                              & - \frac{\mu s}{2} \cdot (k+1)(k+r+1)\cdot \|y_{k+1} - x_{k+1}\|^2                          \nonumber                   \\
                                                              & - \frac{ \mu s}{2} \cdot (2k + r + 3) \cdot \|y_{k+1} - x^{\star}\|^2                         \nonumber                   \\
                                                              & - \frac{s^2(1-Ls) (k+r+1)^2}{2}\left\| \nabla f(y_{k+1}) \right\|^2.                          \label{eqn: iter-mnag1}
\end{align}
To further simplify this expression, we apply the strong convexity inequality~\eqref{eqn: key-inq-smooth}, which gives a lower bound on the gradient norm: 
\begin{equation}
\left\| \nabla f(y_{k+1}) \right\|^2 \geq 2\mu \left( f(x_{k+2}) - f(x^{\star}) \right). 
\label{eqn: key-inq-mnag}
\end{equation}
Substituting this bound~\eqref{eqn: fund-inq-smooth1-new} into~\eqref{eqn: iter-mnag1},  we obtain:
\begin{align}
\mathcal{E}(k+1) - \mathcal{E}(k)  \leq  & - \frac{\mu s(1-Ls)}{4} \cdot s(k+r+1)^2\left( f(x_{k+2}) - f(x^{\star}) \right)          \nonumber                   \\
                                                              & - \frac{\mu s}{2} \cdot (k+1)(k+r+1)\cdot \|y_{k+1} - x_{k+1}\|^2                          \nonumber                   \\
                                                              & - \frac{ \mu s}{2} \cdot (2k + r + 3) \cdot \|y_{k+1} - x^{\star}\|^2                         \nonumber                   \\
                                                              & - \frac{3s^2(1-Ls) (k+r+1)^2}{8}\left\| \nabla f(y_{k+1}) \right\|^2.                        \label{eqn: iter-mnag2}
\end{align}
Next, we estimate the Lyapunov function $\mathcal{E}(k+1)$ defined in~\eqref{eqn: lyapunov-m-nag} using the Cauchy-Schwarz inequality. This gives the upper bound: 
\begin{align}
\mathcal{E}(k+1) \leq  &\;  s(k+2)(k+r+2)\left( f(x_{k+2}) - f(x^{\star}) \right)  + \frac{3(k+1)^2}{2} \|y_{k+1} - x_{k+1}\|^2       \nonumber \\ 
                                    &\;+ \frac{3r^2}{2} \| y_{k+1} - x^{\star} \|^2 + \frac{3s^2(k+r+1)^2}{2} \| \nabla f(y_{k+1}) \|^2.           \label{eqn: lyapunov-estimate-mnag}
\end{align}
By comparing the coefficients on the right-hand sides of inequalities~\eqref{eqn: iter-mnag2} and~\eqref{eqn: lyapunov-estimate-mnag}, we deduce that:
\begin{align}
\mathcal{E}(k+1) - \mathcal{E}(k) & \leq - \mu s \cdot \min\left\{ \frac{1-Ls}{4}, \frac{1}{3}, \frac{2k+r+3}{3r^2}, \frac{1-Ls}{4 \mu s} \right\} \mathcal{E}(k+1)    \nonumber \\
                                                     & \leq - \mu s \cdot \left(\frac{1 - Ls}{4} \right) \cdot \mathcal{E}(k+1),                                                                                     \label{eqn: iter-diff-lyapunov-mnag} 
\end{align}
where the last inequality holds due to the definition of the integer $K$ and the fact that $0 < \mu s < \mu/L \leq 1$. This completes the proof using the Lipschitz gradient condition~\eqref{eqn: grad-lip}. 
\end{proof}

\subsection{Composite optimization via M-FISTA}
\label{subsec: m-fista}
We now extend the linear convergence of \texttt{M-NAG}, as established in~\Cref{thm: m-nag}, to its proximal counterpart,~\texttt{M-FISTA}. As outlined in~\Cref{defn: proximal-value},~\texttt{M-FISTA} utilizes the $s$-proximal value given in~\eqref{eqn: proximal-value} and follows the iteration scheme below, initialized at $y_0 = x_0 \in \mathbb{R}^d$:
\begin{subequations}
\label{eqn: m-fista}
\begin{empheq}[left=\empheqlbrace]{align}
& z_{k}       = P_s(y_{k}),                                                                                                                                                \label{eqn: m-fista-gradient} \\
& x_{k+1}   = \left\{ \begin{aligned} 
                               & z_{k}, && \text{if}\; f(z_{k}) \leq f(x_{k}), \\
                               & x_{k}, && \text{if}\; f(z_{k}) > f(x_{k}),
                               \end{aligned} \right.                                                                                                                          \label{eqn: m-fista-monetone} \\
& y_{k+1}  = x_{k+1} + \frac{k}{k+r+1} (x_{k+1} - x_{k}) + \frac{k+r}{k+r+1} (z_{k} - x_{k+1}),                                        \label{eqn: m-fista-momentum} 
\end{empheq}    
\end{subequations}
where $s>0$ is the step size. To maintain consistency with the \texttt{NAG}-like structure, we introduce an equivalent formulation of~\texttt{M-FISTA} that replaces the $s$-proximal value with the $s$-proximal subgradient, as defined in~\eqref{eqn: subgradient}. This yields the following variant that mirrors the structure of \texttt{M-NAG}:
\begin{subequations}
\label{eqn: m-fista}
\begin{empheq}[left=\empheqlbrace]{align}
& z_{k}       = y_{k} - sG_s(y_{k}),                                                                                                                                   \label{eqn: m-fista1-gradient} \\
& x_{k+1}   = \left\{ \begin{aligned} 
                               & z_{k}, && \text{if}\; f(z_{k}) \leq f(x_{k}), \\
                               & x_{k}, && \text{if}\; f(z_{k}) > f(x_{k}),
                               \end{aligned} \right.                                                                                                                          \label{eqn: m-fista1-monetone} \\
& y_{k+1}  = x_{k+1} + \frac{k}{k+r+1} (x_{k+1} - x_{k}) + \frac{k+r}{k+r+1} (z_{k} - x_{k+1}),                                        \label{eqn: m-fista1-momentum} 
\end{empheq}    
\end{subequations}
where the $s$-proximal subgradient $G_s(y_{k})$ replaces the gradient $\nabla f(y_k)$, thereby enabling extension to composite objectives.

To establish linear convergence for~\texttt{M-FISTA}, we generalize the Lyapunov function~\eqref{eqn: lyapunov-m-nag} to accommodate the proximal structure:
\begin{align}
\mathcal{E}(k) = &\;  s(k+1)(k+r+1)\left( \Phi(x_{k+1}) - \Phi(x^{\star}) \right)                                          \nonumber \\ 
                           &\; +  \frac12\left\| k(y_{k} -x_k) + r(y_k - x^{\star}) - (k + r)s G_s(y_k) \right\|^2.   \label{eqn: lyapunov-m-fista}
\end{align}
This Lyapunov function accounts for the composite nature of the objective $\Phi = f + g$, and the presence of non-smooth components through the proximal subgradient.~\Cref{lem: fund-inq-composite} and~\Cref{lem: key-inq-composite} extend the fundamental inequalities~\eqref{eqn: fund-inq-smooth} and~\eqref{eqn: key-inq-smooth} to the proximal setting, allowing the key inequalities,~\eqref{eqn: fund-inq-xkyk-mnag},~\eqref{eqn: fund-inq-mnag}, and~\eqref{eqn: key-inq-mnag}, used in the analysis of~\texttt{M-NAG} to be generalized to the~\texttt{M-FISTA} setting. With these extension, we can now state the main result:
\begin{theorem}
\label{thm: m-fista}
Let $\Phi = f + g$ with $f \in \mathcal{S}_{\mu,L}^{1}(\mathbb{R}^d)$ and $g \in \mathcal{F}^0(\mathbb{R}^d)$. Given any step size $0 < s < 1/L$, there exists a positive integer $K: = \max\left\{0,  \frac{3r^2 - 4r - 12}{8}\right\}$ such that the iterative sequence $\{x_{k}\}_{k=0}^{\infty}$ generated by~\texttt{M-FISTA}~\eqref{eqn: m-fista},  with any initial $x_0 = y_0 \in \mathbb{R}^d$, satisfies the following inequality: 
\begin{equation}
\label{eqn: m-fista-rate}
\Phi(x_k) - \Phi(x^{\star}) \leq \frac{(r +1)\left( \Phi(x_1) - \Phi(x^{\star})\right) + r^2 L\| x_1 - x^{\star} \|^2}{k(k+r) \left[ 1 + (1 - Ls) \cdot \frac{\mu s}{4} \right]^k}, 
\end{equation}
for any $k \geq \max\left\{ 1, K \right\}$. 
\end{theorem}


\section{Conclusion and future work}
\label{sec: conclusion}

In this paper, we investigate the relationship between the~\texttt{M-NAG} update and the classical \texttt{NAG} scheme through the lens of the implicit-velocity phase-representation. We show that an additional assumption, concerning either the position update or the phase-coupling relation, is necessary, and that each is independently sufficient for~\texttt{M-NAG} to fully recover the~\texttt{NAG} iterates. Furthermore, we demonstrate that the principled approach for constructing a Lyapunov function, as outlined in~\citet{chen2022revisiting}, reveals that the \texttt{M-NAG} update alone is sufficient to establish linear convergence, even without requiring access to the full \texttt{NAG} iteration. To advance this analysis, we introduce a new mixed sequence $\mathbf{S}_k = \mathbf{R}_k - (k+r)s \nabla f(y_k)$, which replaces the earlier sequence $\mathbf{R}_k$ and incorporates forward-indexed gradients for its iterative difference. This leads to a new Lyapunov function that not only guarantees linear convergence for \texttt{NAG} but also omits the kinetic energy term, thereby enabling a direct generalization to \texttt{M-NAG} and overcoming a key technical barrier identified in~\citet{li2024linear}. Finally, by leveraging two proximal inequalities for composite functions developed in~\citet{li2024linear, li2024linear2}, which serves as the proximal analogs of the two classical strong convexity inequalities, we extend the linear convergence guarantees to both~\texttt{FISTA} and~\texttt{M-FISTA}. This extension significantly broadens the applicability of our results and deepens the theoretical understanding of the convergence behavior of these widely used proximal algorithms.

\begin{figure}[htb!]
\centering
\includegraphics[scale=0.24]{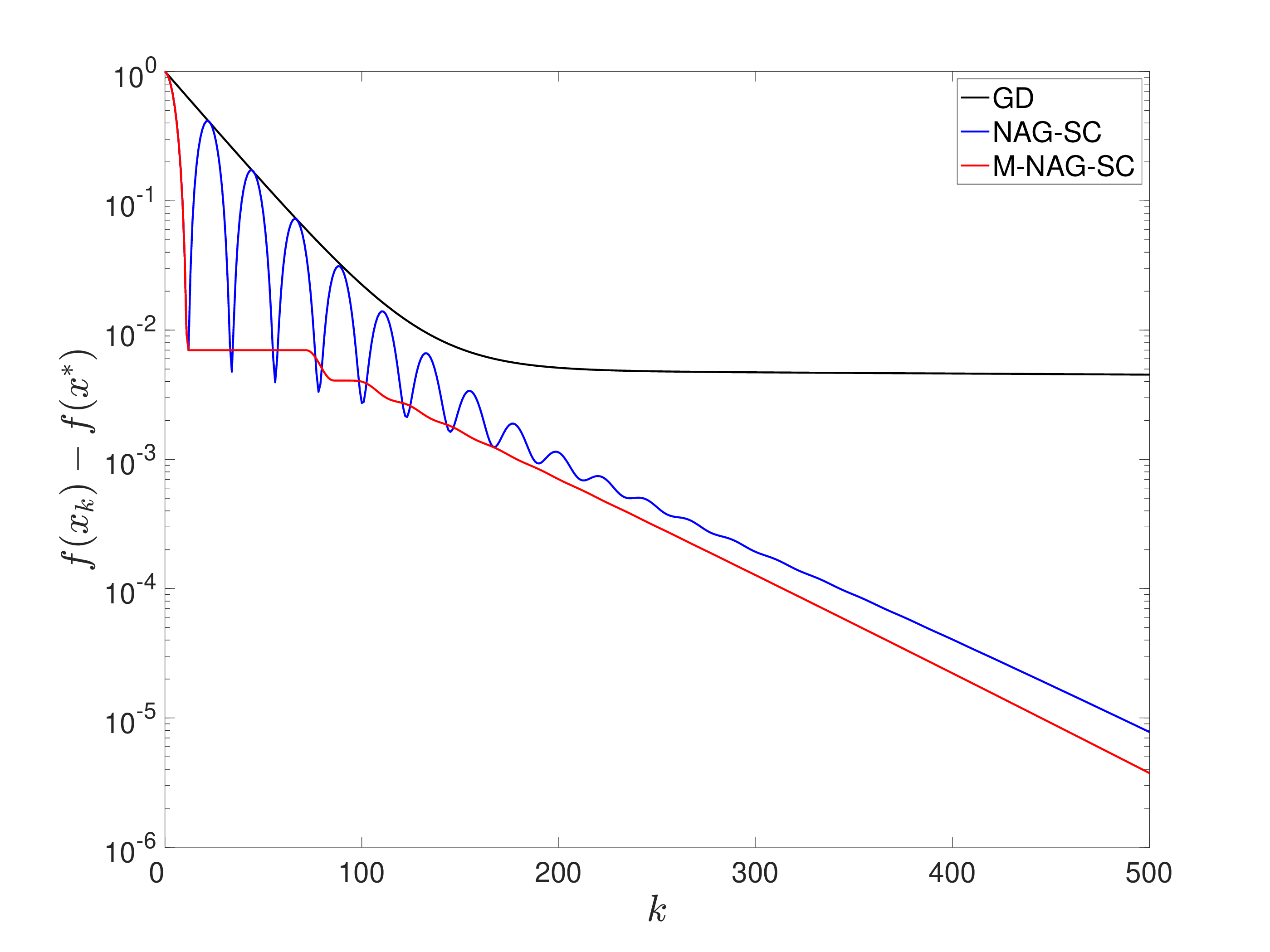}
\caption{Numerical comparison of the progression of function values for vanilla gradient descent,~\texttt{NAG-SC}, and~\texttt{M-NAG-SC}. The step size is set to $s = 0.01$. All optimization algorithms are applied to the quadratic objective function $f(x_1, x_2) = 5 \times 10^{-3}x_1^2 + x_2^2$.}
\label{fig: quadratic_with}
\end{figure}

For $\mu$-strongly convex functions, another variant of Nesterov's accelerated gradient method, referred to as~\texttt{NAG-SC}, was originally proposed by~\citet{nesterov2018lectures}. In~\texttt{NAG-SC}, the momentum coefficient is set to $\frac{1-\sqrt{\mu s}}{1 + \sqrt{\mu s}}$, in contrast to the coefficient $\frac{k}{k+r+1}$ used in the classical~\texttt{NAG}. The iterative scheme for~\texttt{NAG-SC} is:
\begin{subequations}
\label{eqn: nag-sc}
\begin{empheq}[left=\empheqlbrace]{align}
& x_{k+1} = y_{k} - s\nabla f(y_{k}),                                                                                   \label{eqn: nag-sc-gradient}          \\
& y_{k+1} = x_{k+1} + \frac{1-\sqrt{\mu s}}{1 + \sqrt{\mu s}} (x_{k+1} - x_{k}),                 \label{eqn: nag-sc-momentum} 
\end{empheq}    
\end{subequations}
where $s>0$ is the step size. Building on the idea introduced by~\citet{beck2017first}, we propose a monotonic variant of~\texttt{NAG-SC}, denoted as~\texttt{M-NAG-SC}. Its update scheme is given by:
\begin{subequations}
\label{eqn: m-nag-sc}
\begin{empheq}[left=\empheqlbrace]{align}
& z_{k}       = y_{k} - s\nabla f(y_{k}),                                                                                                                                                        \label{eqn: m-nag-sc-gradient} \\
& x_{k+1}   = \left\{ \begin{aligned} 
                               & z_{k}, && \text{if}\; f(z_{k}) \leq f(x_{k}), \\
                               & x_{k}, && \text{if}\; f(z_{k}) > f(x_{k}),
                               \end{aligned} \right.                                                                                                                                                     \label{eqn: m-nag-sc-monetone} \\
& y_{k+1}  = x_{k+1} +\frac{1-\sqrt{\mu s}}{1 + \sqrt{\mu s}}(x_{k+1} - x_{k}) + (z_{k} - x_{k+1}),    \label{eqn: m-nag-sc-momentum} 
\end{empheq}    
\end{subequations}
where $s$ is the step size. We compare the numerical performance of both~\texttt{NAG-SC} and~\texttt{M-NAG-SC} with vanilla gradient descent in~\Cref{fig: quadratic_with}.
An intriguing direction for future research is to investigate whether \texttt{M-NAG-SC}, as described by~\eqref{eqn: m-nag-sc}, can retain the accelerated convergence rate of~\texttt{NAG-SC}. However, this extension poses a nontrivial challenge:  the Lyapunov function developed in~\citet{shi2022understanding}, which incorporates the kinetic energy tern as a critical component, complicates direct generalization to the monotonic setting. Therefore, extending the Lyapunov analysis to cover~\texttt{M-NAG-SC} remains an open and promising challenge.


\section*{Acknowledgements}
We thank Bowen Li for his helpful and insightful discussions. Mingwei Fu acknowledges partial support from the Loo-Keng Hua Scholarship of CAS. This project was partially supported by the startup fund from SIMIS and Grant No.12241105 from NSFC.

\bibliographystyle{abbrvnat}
\bibliography{sigproc}

\end{document}